\documentclass[11pt]{article}
\input epsf.tex
\usepackage{latexsym}
\usepackage{amsmath,amsthm,amsfonts,amssymb}
\usepackage{epsfig}

\newtheorem{pro}{Proposition}[section]
 \newtheorem{thm}[pro]{Theorem}
 \newtheorem{lem}[pro]{Lemma}

 \newtheorem{cor}[pro]{Corollary}
\def\A{{\cal A}}
 
 \def\C{{\Bbb C}}
\def\coupling{{\cal C}}
\def\E{{\Bbb E}}
 \def\H{{\Bbb H}}
 \def\G{{\cal G}}

 \def\chix{{\raise.5ex\hbox{$\chi$}}}

\begin{document}
\title{Couplings of Uniform Spanning Forests}
\author{Lewis Bowen \footnote{Research supported in part by NSF Vigre Grant No. DMS-0135345} }

\maketitle

\begin{abstract}
We prove the existence of an automorphism-invariant coupling for the wired and the free uniform spanning forests on connected graphs with residually amenable automorphism groups.
\end{abstract}



\section{Introduction}
All the graphs in this paper will be assumed to be locally finite with at most a countable number of vertices. Given a (connected) graph $\Gamma$, there exists two natural random subgraphs on $\Gamma$ known as the Free Uniform Spanning Forest (FSF) and the Wired Uniform Spanning Forest (WSF) (see [BLPS] or [L1] for introductory material and more). They are defined as follows. Consider any increasing sequences $\Gamma_i$ of finite connected subgraphs of $\Gamma$ whose union is all of $\Gamma$. Then, as shown in [Pe], the uniform distribution on spanning trees on $\Gamma_i$ converges and is, by definition the distribution of FSF. To obtain WSF, identify all of the boundary vertices of $\Gamma_i$ to obtain a graph $\tilde{\Gamma_i}$. As proven in [Pe], the uniform distribution on spanning trees on $\tilde{\Gamma_i}$ converges and the limiting distribution is the distribution of WSF (by definition). 

In [BLPS], it is shown that there exists a monotone coupling $\mu$ of WSF and FSF (this follows from [FM]), i.e. $\mu$ is a probability measure on $2^{E(\Gamma)}\times 2^{E(\Gamma)}$ (where $E(\Gamma)$ denotes the edge set of $\Gamma$) whose projection onto the first factor is the distribution of WSF, whose projection onto the second factor is the distribution of FSF, and such that $\mu$ is concentrated on the set of all pair $(T_1,T_2) \in 2^{E(\Gamma)} \times 2^{E(\Gamma)}$ such that $T_1 \subset T_2$. The proof, however, is purely an existence proof. The authors of [BLPS] posed the following problem: does there exist a natural or explicit monotone coupling of WSF and FSF? In particular, does there exist an automorphism invariant coupling of WSF and FSF?

We show if $\Gamma$ is a (connected, locally finite) graph and $G < Aut(\Gamma)$ is a residually amenable group of automorphisms of $\Gamma$ then there exists a coupling of WSF and FSF that is invariant under elements of $G$. We note that all finitely generated linear groups are residually amenable. See section 2 for more details on residual amenability.

We mention that for many graphs $\Gamma$, WSF = FSF. This is true if and only if there are no nonconstant harmonic Dirichlet functions on $\Gamma$ (see [BLPS]). In particular, WSF=FSF whenever $\Gamma$ is the Cayley graph of a cocompact lattice $G$ in $\E^n$ (for all $n$) or in $\H^n$ (for $n > 2$). However, WSF $\ne$ FSF when $G$ is a cocompact Fuchsian group.

\section{Background on residually amenable groups}

We will use the following definition of amenable (see [Z]):

A topological group $G$ is {\bf amenable} if and only if for every compact metrizable space $Z$ such that $G$ acts on $Z$ by homeomorphisms and the map $G \to Homeo(Z)$ is continuous, there exists a Borel probability measure $\mu$ on $Z$ that is invariant under the action of $G$. By invariant, we mean that for all Borel sets $E \subset Z$ and for all $g \in G$, $\mu(gE)=\mu(E)$. The topology on $Homeo(Z)$ is the compact open topology. We note that finite groups and finitely generated abelian groups are amenable (given the discrete topology).

If $P$ is a property of groups, then a topological group $G$ is said to be {\bf residually $P$} if for every $g \in G-\{1\}$, there exists a surjective continuous homomorphism $\phi: G \to H$ such that $H$ has property $P$ and $\phi(g)$ is not the identity. 

We need the following two results from [Z] Prop. 4.1.6. First, if $G$ is amenable and $H$ is a closed subgroup of $G$, then $H$ is amenable. Second, if $K$ is a closed normal subgroup of $G$ then $G$ is amenable if and only if both $K$ and $G/K$ are amenable. Now suppose that $K_1$ and $K_2$ are closed normal subgroups of $G$ such that both $G/K_1$ and $G/K_2$ are amenable. We claim that $G/(K_1 \cap K_2)$ is amenable. To see this, note that $K_1/(K_1 \cap K_2)$ is a closed normal subgroup of $G/(K_1 \cap K_2)$ and $[G/(K_1 \cap K_2)]/[K_1/(K_1 \cap K_2)]$ is isomorphic to $G/K_1$ which is amenable. Also $K_1 /(K_1 \cap K_2)$ maps injectively into $G/K_2$ (under the quotient map) so it is amenable too (by the first result). The second result now implies $G/(K_1 \cap K_2)$ is itself amenable. It follows easily that if $G$ is separable (i.e. $G$ contains a countable dense subset) then $G$ is residually amenable if and only if there exists a decreasing sequence of closed normal subgroups $G_i$ in $G$ such that $G/G_i$ is amenable and $\cap_i \, G_i = \{1\}$. The next lemma is presumably well-known but we did not find it in the literature.

\begin{lem}
If $\Gamma$ is a locally finite graph (with at most a countable number of vertices) and $G$ is a group of automorphisms of $\Gamma$ then $G$ is separable (in the open-compact topology).
\end{lem}

\begin{proof}
Let $V$ denote the vertex set of $\Gamma$. We let $B_r(v)$ denote the radius $r$ neighborhood of a vertex $v \in V$. Since $\Gamma$ is locally finite, for every integer radius $r>0$ and $(v,w) \in V \times V$ there are only finitely many bijections from $B_r(v)$ to $B_r(w)$. So there exists a countable set $S \subset G$ such that for any $(v,w) \in V \times V$ if there is a $g \in G$ such that $g(v)=w$ then there exists a $g_r \in S$ that agrees with $g$ on the radius $r$ neighborhood of $v$. In the open-compact topology $g_r$ converges to $g$ as $r \to \infty$. So $S$ is dense in $G$.
\end{proof}

We let $Aut(\Gamma)$ denote the full automorphism group of a graph $\Gamma$.

\begin{cor}
If $G<Aut(\Gamma)$ is residually amenable then there exists a decreasing sequence $\{G_i\}$ of closed normal subgroups of $G$ such that $G/G_i$ is amenable and $\cap_i G_i=\{1\}$.
\end{cor}

We note that since finite groups are amenable, any residually finite group is residually amenable. A well-known result due to Mal'cev [Ma] (or see theorem 2.7 [We]) states that if $G$ is a finitely generated subgroup of $GL_n(\mathbb{C})$ then $G$ is residually finite. In particular, all discrete groups of isometries of hyperbolic space are residually finite. It is an open problem whether or not all word hyperbolic groups are residually finite (see [KW]). Olshanskii has proved [O] that in a certain sense, almost every finitely generated group is word hyperbolic. So, if all word hyperbolic groups are residually finite then almost every (finitely generated) group is residually amenable. 

There is an example, due to E.A. Scott [S] of a finitely generated infinite simple group that has a free nonabelian subgroup. Such a group cannot be residually amenable. There are groups that are residually amenable but not residually finite. For example, the Baumslag-Solitar group $BS(2,3)= <x,y \, | \, xy^2x^{-1}=y^3>$ is one.

\section{Background on Determinantal Probability Measures}

Let $E$ be any set (with at most countable cardinality), and $S < l^2(E)$ be a closed subspace of $l^2$-summable functions on $E$. (In general, we write $S<T$ to mean $S$ is a closed subspace of $T$). There exists a natural probability measure ${\bf P}^S$ on $2^E$ defined by 
\begin{equation}
{\bf P}^S(B \subset \mathfrak{B}) = Det[<P_S(\chi^e),\chi^f>_{e,f \in B}].
\end{equation}

${\bf P}^S$ is said to be a determinantal probability measure (see [L2] for the history and background of such measures). Here, $B$ is a finite subset of $E$ and ${\bf P}^S(B \subset \mathfrak{B})$ is by definition the probability with respect to ${\bf P}^S$ that $B$ is contained in a random subset of $\Gamma$, $P_S$ denotes projection onto $S$, $\chi^e$ denotes the characteristic function of an element $e \in B$. The right hand side is the determinant of the matrix with columns and rows indexed by elements of $B$ whose $(e,f)$-entry is equal to the $l^2$-inner product $<P_S(\chi^e),\chi^f>$. Using the inclusion-exclusion principle, it can be shown that this determines a probability measure. A determinantal probability measure is defined as one that arises from this construction (see [L2] for introductory material and more). We let $\G_S$ denote a random subset with distribution ${\bf P}^S$.

It is shown in [BLPS] that the distributions of WSF and FSF are determinantal. To be precise, let $\Gamma$ be a graph. We let $l^2_-(\Gamma)$ be the space of $l^2$-summable antisymmetric functions $f$ on the (directed) edge set of $\Gamma$. By antisymmetric we mean that if $e$ is an edge of $\Gamma$ and $\check{e}$ is equal to $e$ with the opposite orientation, then $f(e)=-f(\check{e})$. For any edge $e$ of $\Gamma$, let ${\bf 1}_e$ be the characteristic function of $e$. We let $\chi^e = {\bf 1}_e - {\bf 1}_{\check{e}}$ be the unit flow along the (directed) edge $e$.

For any vertex $v$ of $\Gamma$, define the star of $v$ is defined to be $\Sigma_e \, \chi^e$ where the sum is over all directed edges $e$ with initial endpoint equal to $v$. We let $\bigstar$ be the space in $l^2_-(\Gamma)$ generated by all stars of all vertices.

We say that $(e_1,..,e_k)$ is a cycle of edges in $\Gamma$, if the destination vertex of $e_i$ equals to the initial vertex of $e_{i+1}$ for all $i$ mod $k$. In this case, we say that the function $\Sigma_i \, \chi^e \in l^2_-(\Gamma)$ is a cycle. We let $\diamondsuit$ be the subspace of $l^2_-(\Gamma)$ generated by all cycles in $\Gamma$. It is shown in [BLPS] that

\begin{equation}
l^2_-(\Gamma)= \bigstar \oplus \diamondsuit \oplus \nabla HD
\end{equation}

\noindent where $\nabla HD$ is the gradient of the space of all harmonic Dirichlet functions on $\Gamma$. It is proven in [BLPS] that the distributions of WSF and FSF are the determinantal probability measures associated to $\bigstar$ and $\bigstar \oplus \nabla HD = \diamondsuit^\perp$ respectively. To be precise, if $\vec{B}$ is a set of directed edges such that if $e \in \vec{B}$ then $\check{e} \notin \vec{B}$ then let $B$ denote the same set of edges of $\vec{B}$ but without orientation. Then we define ${\bf P}^S$ by

\begin{equation}
{\bf P}^S(B \subset \mathfrak{B}) = Det[<P_S(\chi^e), \chi^f>_{e,f \in \vec{B}}].
\end{equation}

\noindent where $S$ is a subspace of $l^2_-(\Gamma)$. It follows from elementary properties of deteminants that this is well-defined (i.e. it is independent of the choice of $\vec{B}$). It is shown in [BLPS] that ${\bf P}^\bigstar$ is the distribution of WSF and ${\bf P}^{\diamondsuit^\perp}$ is the distribution of FSF. For our purposes, it suffices to assume this as the definition of the distributions of WSF and FSF. One easily deduces that WSF=FSF if and only if $\nabla HD = 0$. This occurs, for example when $\Gamma$ is the Cayley graph of an amenable group or of a Kazhdan group. See Bekka and Valette [BV] (Theorem D) for a more complete list of groups whose Cayley graphs satisfy $\nabla HD = 0$.

We will need the following facts from [L2] (Theorem 5.2). If $S < T < l^2_-(\Gamma)$ and $S$ and $T$ are closed then for any increasing event $\A \subset 2^{E(\Gamma)}$ ($\A$ is increasing whenever $A \in A$ and $A \subset B \subset E(\Gamma)$ implies $B \in \A$)
 ${\bf P}^{S}(\A) \le {\bf P}^T(\A)$. By Strassen's theorem [St], this implies the existence of a monotone coupling between ${\bf P}^S$ and ${\bf P}^T$.

\section{the theorem}

In this section we state and prove the main theorem.

\begin{thm}
If $\Gamma$ is a connected locally finite graph (with at most a countable number of vertices) and $G < Aut(\Gamma)$ is residually amenable (with respect to the compact-open topology) then there exists a $G$-invariant monotone coupling between the WSF and the FSF on $\Gamma$.
\end{thm}

Since $G$ is residually amenable and separable, there exists a decreasing sequence $\{G_i\}$ of closed normal subgroups of $G$ such that $G/G_i$ is amenable for all $i$ and $\cap_i G_i = \{1\}$. 

We let $\Gamma/G_i$ be the quotient graph. Its vertices (edges) are equivalence classes of vertices (edges) in $\Gamma$ where $v$ is equivalent to $w$ if there is a $g \in G_i$ such that $gv=w$. If $[v]$ and $[w]$ are vertices of $\Gamma/G_i$, then there is an edge between them for every equivalence class $[e]$ of edges such that one endpoint of $e$ is in $[v]$ and the other is in $[w]$. Consider $l^2_-(\Gamma/G_i)$. We let $\bigstar_i$ be the subspace generated by all the stars of $\Gamma/G_i$ and $\diamondsuit_i$ be the subspace generated by all the cycles of $\Gamma/G_i$. We let $\nabla HD_i$ denote the gradient of the Harmonic Dirichlet functions on $\Gamma/G_i$. So,

\begin{equation}
l^2_-(\Gamma/G_i) = \bigstar_i \oplus \diamondsuit_i \oplus \nabla HD_i.
\end{equation}

Let $\pi_i : \Gamma \to \Gamma/G_i$ denote the quotient map. Suppose $c$ is a cycle in $l^2_-(\Gamma/G_i)$. Then $c = \Sigma_{i=1}^n \chi^{e_i}$ where $e_1,e_2,..,e_n$ is a sequence of directed edges such that the terminal vertex of $e_i$ is equal to the initial vertex of $e_{i+1}$ (mod $n$). We will say that $c$ is a true cycle if every directed component of $\pi_i^{-1}(\{e_1,..,e_n\})$ forms a (finite) cycle in $\Gamma$. We let $C_i$ denote the closed subspace of $\diamondsuit_i$ spanned by all the true cycles and let $H_i$ denote its orthocomplement in $\diamondsuit_i$. Now we have

\begin{equation}
l^2_-(\Gamma/G_i) = \bigstar_i \oplus C_i \oplus H_i \oplus \nabla HD_i.
\end{equation}

Now let $W_i$ be the random subgraph of $\Gamma/G_i$ associated to $\bigstar_i$. Let $F_i$ be the random subgraph associated to $\bigstar_i \oplus H_i \oplus \nabla HD_i = C_i^\perp$.

 The strategy of the proof is as follows. We show that in a certain sense, a lift of $W_i$ converges to WSF and a lift of $F_i$ converges to FSF. Using the amenability of $G/G_i$, we can average any coupling of $W_i$ and $F_i$ over $G/G_i$ to obtain a $G/G_i$-invariant coupling. By a compactness argument, a limit point of the sequence of lifted couplings exists and we will show that any limit point is a $G$-invariant monotone coupling of WSF and FSF.

We need to introduce some topological spaces. First, we let $Z = 2^{E(\Gamma)}$ be the space of subgraphs of $\Gamma$. It is compact in the product topology. The group $G$ acts on $Z$ in the natural way and this action is bi-continuous. We let $M_Z$ denote the set of Borel probability measures on $Z$ and $M^i_Z$ denote the subset of $M_Z$ of measures that are $G$-invariant. By invariance we mean if $\mu \in M^i_Z$, $E \subset Z$ and $g \in G$ then $\mu(gE) =\mu(E)$. We say that a sequence $\{\mu_i\}$ in $M_Z$ converges to $\mu$ in the weak* topology if

\begin{equation}
\int_Z f \, d\mu_i \to \int_Z f \, d\mu
\end{equation}

\noindent for every continuous function $f: Z \to \C$. It follows from standard functional analysis that both $M_Z$ and $M^i_Z$ are compact under the weak* topology. Note that for every finite subset $B \subset E(\Gamma)$, the function $f_B: Z \to \C$ defined by $f_B(z)=1$ if $B \subset z$ and $f_B(z)=0$ otherwise is continuous. It is not hard to show that if

\begin{equation}
\int_Z f_B \, d\mu_i \to \int_Z f_B \, d\mu
\end{equation}

\noindent for every finite set $B \subset E(\Gamma)$ then $\mu_i$ converges to $\mu$ in the weak* topology. In this case, if $\G_i$ is a random subgraph with distribution $\mu_i$ and $\G$ is a random subgraph with distribution $\mu$ then we will say that $\G_i$ converges to $\G$ (weak*). Similarly, we define $M^i_{Z \times Z}$ to be the space of $G$-invariant Borel probablity measures on $Z \times Z$. It is compact under the weak* topology.

If $\{S_i\}$ is a sequence of closed subspaces in $l^2_-(\Gamma)$, we will say that $P_{S_i} \to P_S$ in the strong operator topology (SOT) if for every $f \in l^2_-(\Gamma)$, $||P_{S_i}(f) - P_S(f)|| \to 0$. We sometimes express this by writing $S_i \to S$ (SOT). It can be shown that the set of subspaces of $l^2_-(\Gamma)$ is compact under the SOT. Note that if $P_{S_i} \to P_S$ (SOT) then ${\bf P}^{S_i} \to {\bf P}^S$ in the weak* topology. It is well known that if $S_i < l^2_-(\Gamma)$ are closed subspaces and $S_i \nearrow S$ (meaning that $S_i \subset S_{i+1}$ and $\cup S_i$ is dense in $S$) then $P_{S_i}(T) \to P_{S}(T)$ (SOT) for any closed subspace $T$ (in $l^2_-(\Gamma)$).

 If $S < l^2_-(\Gamma/G_i)$, let ${\bf \tilde{P}}^S$ be the distribution of the random subgraph $\pi_i^{-1}(\G_S)$. For ease of notation, we also denote $\pi_i^{-1}(W_i)$ by $\tilde{W_i}$ and $\pi_i^{-1}(F_i)$ by $\tilde{F_i}$. We will show that $\tilde{W_i} \to$ WSF and $\tilde{F_i} \to$ FSF. In order to do this, we introduce sequences of random subgraphs of $\Gamma$ that equal $\tilde{W_i}$ or $\tilde{F_i}$ on fundamental domains of $G_i$ and have the advantage that their distributions are determinantal. For this we let $D_i$ be a connected subgraph of $\Gamma$ such that the covering map $\pi_i$ restricted to the edge set of $D_i$ is bijective. Also assume that $D_i \subset D_{i+1}$ and that $\cup_i D_i = \Gamma$. This is possible since $\cap_i G_i = \{1\}$.

If $\G_S$ is a random subgraph of $\Gamma/G_i$, let $\hat{\G_S}$ be the random subgraph of $\Gamma$ that is contained in $D_i$ and is equal to $\tilde{\G_S}$ on $D_i$. Note that $\pi_i$ restricted to $D_i$ induces an isomorphism $\pi_{i*}: l^2_-(D_i) \to l^2_-(\Gamma/G_i)$. If $S < l^2_-(\Gamma/G_i)$ is a (closed) subspace, let ${\hat S}$ denote $\pi_{i*}^{-1}(S)$. It is easy to check that ${\hat \G_S}$ is the random subgraph associated to ${\hat S}$ (i.e. ${\hat \G_S} = \G_{\hat S}$).

\begin{lem}
Suppose for each $i$, $\G_i$ is a random subgraph on $\Gamma/G_i$. Then ${\tilde \G_i} \to \G$ (weak*) if and only if ${\hat \G_i} \to \G$ (weak*).
\end{lem}

\begin{proof}
Let $B$ be a finite set of edges of $\Gamma$. Suppose that $\tilde{\G_i} \to \G$ (weak*). Then, by definition of weak* convergence,

\begin{equation}
\int_Z f_B \, d{\tilde \mu}_i \to \int_Z f_B \, d\mu
\end{equation}

\noindent where ${\tilde \mu}_i$ is the distribution of $\tilde{\G_i}$ and $\mu$ is the distribution of $\G$. There exists an $N$ such that $i>N$ implies that $B \subset D_i$. Since $\hat{\G_i}$ is equal to $\tilde{G_i}$ on $D_i$, it follows that  

\begin{equation}
\int_Z f_B \, d{\hat \mu}_i \to \int_Z f_B \, d\mu
\end{equation}

\noindent where ${\hat \mu}_i$ is the distribution of ${\hat \G_i}$. Since $B$ is arbitrary, ${\hat \G_i} \to \G$ (weak*). The proof in the other direction is similar. 

\end{proof}

\begin{thm}
${\tilde W_i} \to $ WSF and ${\tilde F_i} \to $ FSF (weak*).
\end{thm}

\begin{proof}

By the above lemma, it suffices to prove that ${\hat W_j} \to $ WSF and ${\hat F_j}$ to FSF. To this end, it suffices to prove that ${\hat \bigstar_j} \to \bigstar$ and $\widehat{C_j^\perp} \to \diamondsuit^\perp$. Let $B$ be a finite set of (directed) edges of $\Gamma$ and let $S$ be the subspace generated by $\chi^e$ for $e \in B$. 

We show first that $P_S(\bigstar) < P_S({\hat \bigstar_j})$ for all large enough $j$. Suppose that $v$ is a vertex incident to an edge $e \in B$. Let $f_v$ be the star of $v$. It suffices to show that $P_S(f_v) \in P_S({\hat \bigstar_j})$ for all large enough $j$. So let $j$ be large enough so that all edges incident to $v$ are in $D_j$. Let $g_v$ be the star of $\pi_j(v)$ in $\Gamma/G_j$. By definition $g_v \in \bigstar_j$. It is clear that $f_v = \pi_{j*}^{-1}(g_v)$. So  $P_S(f_v) = P_S(\pi_{j*}^{-1}(g_v))$ which implies that $P_S(f_v) \in P_S({\hat \bigstar_j})$. Since $v$ is arbitrary, the claim is proven. 

To show that $P_S({\hat \bigstar_j}) < P_S(\bigstar)$ for $j$ large enough we follow the same argument as above in reverse. So, $P_S({\hat \bigstar_j}) = P_S(\bigstar)$ for all $j$ large enough. Similarly, it can be shown that  $P_S({\hat C_j}) = P_S(\diamondsuit)$ for all $j$ large enough.

 We let $S_i$ be the subspace generated by the functions $\chi^e$ for $e \in D_i$. Then $S_i \nearrow l^2_-(\Gamma)$ so $P_{S_i}(\bigstar) \to \bigstar$ (SOT). Since $P_{S_i}(\bigstar)=P_{S_i}({\hat \bigstar_j})$ for $j$ large enough, we have that for any subsequential limit $T$ for $\{{\hat \bigstar_j}\}$, $P_{S_i}(T) \to \bigstar$ (SOT). But $P_{S_i}(T) \to T$ (SOT) so $T=\bigstar$, i.e. ${\hat \bigstar_j} \to \bigstar$ (SOT). This implies that ${\hat W_i} \to$ WSF (weak*).

Similarly, it can be shown that ${\hat C_i} \to \diamondsuit$ (SOT). This implies that $({\hat C_i})^\perp \to \diamondsuit^\perp$ (SOT).  Since $S_i \nearrow l^2_-(\Gamma)$, it follows that $P_{S_i}(({\hat C_i})^\perp) \to \diamondsuit^\perp$. But note that $P_{S_i}(({\hat C_i})^\perp) = \widehat{C_i^\perp}$. Hence we have shown that $\widehat{C_i^\perp} \to \diamondsuit^\perp$. This implies that ${\hat F_i} \to$ FSF (weak*) since $\widehat{C_i^\perp}$ is the distribution of ${\hat F_i}$. Now we are done.

\end{proof}

\begin{proof}( of theorem 4.1)

By Strassen's lemma and [L2] Theorem 6.1, for all $i$, there exists a monotone coupling of $W_i$ and $F_i$. 
 The space $\coupling_i$ of all such monotone couplings is compact under the weak* topology. $G/G_i$ acts on this space by $gG_i\mu(F) \to \mu(gG_iF)$ for $\mu \in \coupling_i$ and $F \subset 2^{E(\Gamma/G_i)} \times 2^{E(\Gamma/G_i)}$. Since $G/G_i$ is amenable there exists a $G/G_i$-invariant measure $\sigma_i$ on $\coupling_i$. Let $\mu_i$ be defined by 

\begin{equation}
\mu_i(F) = \int_{\coupling_i} \nu(F) \, d\sigma_i(d\nu)
\end{equation}
 for any Borel set $F \subset 2^{E(\Gamma/G_i)}$. Because $\sigma_i$ is $G/G_i$-invariant, $\mu_i$ is a $G/G_i$-invariant monotone coupling of $W_i$ and $F_i$.

Since $\pi_i: \Gamma \to \Gamma/G_i$ is equivariant with respect to the $G$-action, and $\mu_i$ is invariant under $G/G_i$, we  may pull $\mu_i$ back via $\pi_i$ to a $G$-invariant monotone coupling ${\tilde \mu_i}$ between ${\tilde W_i}$ and ${\tilde F_i}$.  

 The space of $G$-invariant monotone couplings is naturally a closed subspace of $M^i_{Z \times Z}$. This space is compact, so there exists a weak* limit point $\mu$ of ${\tilde \mu}_i$.  Since $\mu_i \to \mu$, the first and second marginals of $\mu_i$ must converge to the first and second marginals of $\mu$ respectively. Thus $\mu$ is a coupling between WSF and FSF.

\end{proof}

\section{speculations}

\begin{enumerate}

\item Suppose that $\Gamma$ is a planar graph embedded in the hyperbolic plane so that $Aut(\Gamma)$ acts by hyperbolic isometries. We may then take $G = Aut(\Gamma)$ and choose $G_i$ so that $G/G_i$ is finite. It can be shown that the distribution of $F_i$ is the uniform distribution on subgraphs $\G$ of $\Gamma/G_i$ such that $\G$ contains a spanning tree and $\H^2/G_i - \G$ is a topological disk. If there exists a ``natural'' or explicit coupling between $W_i$ and $F_i$ then one may hope to obtain a natural coupling between WSF and FSF as a limit of lifts of couplings of $W_i$ and $F_i$. It seems that it should be possible to find such a coupling (if it exists) if only because $\Gamma/G_i$ is a finite graph.

\item Suppose that $S$ is a $G$-invariant subspace of $l^2_-(\Gamma)$. Does there exist $G/G_i$-invariant subspaces $S_i$ in $l^2_-(\Gamma/G_i)$ such that ${\hat S_i} \to S$?

\end{enumerate}

{\bf Acknowledgements}
I am grateful to Henry Cohn, Scott Sheffield, Russell Lyons and Yuval Peres for valuable conversations.

Mathematics Department
University of California
Davis, CA
95616,
lbowen@math.ucdavis.edu


\begin{thebibliography}{xxx}


 \bibitem[BLPS]{BLPS} I. Benjamini, R. Lyons, Y. Peres and O. Schramm, \textit{Uniform Spanning Forests}, Ann. Probab. 29 (2001), no. 1, 1--65.

\bibitem[BV]{BV} Bekka, M. E. B., Valette A. \textit{Group cohomology, harmonic functions and the first $L^2$-Betti number}, Potential Anal. (1997), 6, 313-326.
\bibitem[FM]{FM}Feder, T. and Mihail, M. (1992). \textit{Balanced Matroids}, Proc. 24th Annual ACM Sympos. Theory Computing (Victoria, BC, Canada), pp 26-38. ACM Press, New York.

\bibitem[L1]{L1} Lyons, R. \textit{A bird's-eye view of uniform spanning trees and forests}, in Microsurveys in Discrete
   Probability, D. Aldous and J. Propp (eds.), Amer. Math. Soc., Providence, RI, 1998, pp.
   135--162.

 \bibitem[L2]{L2} Lyons, R. \textit{Determinantal Probability Measures}, preprint.

 

\bibitem[KW]{KW} Kapovich, Ilya; Wise, Daniel T. \textit{The equivalence of some residual properties of word-hyperbolic groups.} J. Algebra 223 (2000), no. 2, 562--583.

\bibitem[Ma]{Ma} Mal'cev, A.I. \textit{On faithful representations of infinite groups of matrices} (Russian). Mat Sb. 8, 405-422 (1940). = Amer. Math. Soc. Transl. (2) 45, 1-18 (1965).

\bibitem[O]{O} Olshanskii, A. Yu. \textit{Almost every group is hyperbolic.} Intern. J. Alg. Comp. 2 (1992), no.1, 1-17.

\bibitem[Pe]{Pe} Pemantle, Robert. \textit{Choosing a spanning tree for the integer lattice uniformly.} Ann. Probab. 19 (1991), 1559-1574.

\bibitem[Pi]{Pi} Pier, Jean-Paul. \textit{Amenable locally compact groups.} Pure and Applied Mathematics. A Wiley-Interscience Publication. 
John Wiley \& Sons, Inc., New York, 1984. x+418 pp.

\bibitem[S]{S} Scott, E.A. \textit{A tour around finitely presented infinite simple groups}, Algorithms and classification in combinatorial group theory (Berkeley, 1989), pp. 83-119, Springer-Verlag, New York, 1992.

\bibitem[St]{St} Strassen, V. \textit{The existence of probability measures with given marginals.} Ann. Math. Statist. 36, (1965), 423-439.

\bibitem[We]{We} Wehrfritz, B.A.F. \textit{Infinite Linear Groups}. Springer-Verlag, New York, Heidelberg, Berlin, 1973. 

\bibitem[Z]{Z} Zimmer, Robert J. \textit{Ergodic theory and semisimple groups.} Monographs in Mathematics, 81. Birkh\"auser Verlag, Basel, 1984. x+209 pp.

 \end{thebibliography}
\end{document}